%% file: main.tex
\documentclass[10pt,letterpaper]{article}
\usepackage[top=0.85in,left=1.5in,right=1.5in,footskip=0.75in,marginparwidth=1.5in]{geometry}

\usepackage[utf8]{inputenc}
\usepackage{epsfig} 
\usepackage{graphicx}
\usepackage{subfigure}
\usepackage{mathtools}
\usepackage{algpseudocode}
\usepackage{algorithm}
\usepackage{amsmath}
\usepackage{bm}
\usepackage{amsfonts}
\usepackage{color}
\usepackage{soul}
\usepackage{xcolor}
\usepackage{verbatim}
\usepackage{cite}
\usepackage{mathabx}
\usepackage{nameref,hyperref}

\usepackage[right]{lineno}

\usepackage{microtype}
\DisableLigatures[f]{encoding = *, family = * }

\usepackage{changepage}

\usepackage[aboveskip=1pt,labelfont=bf,labelsep=period,singlelinecheck=off]{caption}

\makeatletter
\renewcommand{\@biblabel}[1]{\quad#1.}
\makeatother

\usepackage{lastpage,fancyhdr,graphicx}
\usepackage{epstopdf}
\fancyheadoffset[L]{2.25in}
\fancyfootoffset[L]{2.25in}

\newcommand{\bx}{\mathbf{x}}

\newcommand{\by}{\mathbf{y}}

\newcommand{\bxi}{\boldsymbol{\xi}}

\newcommand{\E}{\mathbb{E}}
\newcommand{\R}{\mathbb{R}}
\newcommand{\GP}{\operatorname{GP}}

\newcommand{\calN}{\mathcal{N}}

\newcommand{\sref}[1]{Sec.~\ref{sec:#1}}
\newcommand{\qref}[1]{Eq.~(\ref{eqn:#1})}

\newcommand{\fref}[1]{Fig.~\ref{fig:#1}}
\newcommand{\fcite}[1]{\cite{#1}}

\usepackage{color}

\definecolor{Gray}{gray}{.25}

\usepackage{wrapfig}
\reversemarginpar

\begin{document}
\vspace*{0.1in}

\begin{center}
{\Large
\textbf\newline{Stochastic Multi-objective Optimization on  a  Budget:  
Application  to  multi-pass  wire drawing with quantified uncertainties}
}
\newline
\\
Piyush Pandita\textsuperscript{1},
Ilias Bilionis\textsuperscript{1,*},
Jitesh Panchal\textsuperscript{1},
B.P. Gautham\textsuperscript{2},
Amol Joshi\textsuperscript{2},
Pramod Zagade\textsuperscript{2},
\\
\bigskip
{1} School of Mechanical Engineering, Purdue University, West Lafayette, Indiana 47907
\\
{2} Tata Research Development and Design Centre, Tata Consultancy Services, Pune,
India
\\
\bigskip
* ibilion@purdue.edu

\end{center}

\input{abs}

\input{intro}

\section{Methodology}
\label{sec:metho}

Let $X$ denote the set of feasible designs and $(\Omega,\mathcal{F},\mathbb{P})$ be a probability space.
We assume that $X$ is a closed and bounded set of a Euclidean space.
We have $m$ stochastic quantities of interest (QoIs) which we represent as Borel-measurable functions $o_i:X\times\Omega\rightarrow\mathbb{R}, i=1,\dots,m$.
Our goal is to find designs $\bx\in X$ that maximize the expectations of these QoIs over $\omega\in\Omega$, i.e., we wish to maximize $O_i(\bx):=\mathbb{E}[o_i(\bx,\omega)] := \int o_i(\bx,\omega)d\mathbb{P}(\omega)$.
We say that $\bx\in X$ \emph{dominates} $\bx'\in X$, and write $\bx\succcurlyeq\bx'$, if and only if
\begin{equation}
  O_i(\bx) \ge O_i(\bx'),\forall i=1,\dots,m.
\end{equation}
We say that $\bx$ \emph{strictly dominates} $\bx'$, and write $\bx\succ\bx'$ if and only if $\bx\succcurlyeq\bx'$ and there exists $i\in\{1,\dots,m\}$ such that $\mathbb{E}[o_i(\bx,\omega)] > \mathbb{E}[o_i(\bx',\omega)]$.

We wish to characterize the \emph{set of optimal designs}, otherwise known as the \emph{Pareto-efficient frontier}, induced by the preference relation `$\succcurlyeq$'.
In words, the Pareto-efficient frontier, $P_O$, is the set of achievable objectives that are not dominated.
Since $P_O$ has Lebesgue measure zero, working with it directly is problematic.
Instead, we will work with the \emph{attained set}, $A_O$, which is defined as the set of achievable objectives that are strictly dominated.
$P_O$ is simply part of the boundary of $A_O$.

We now proceed to the exact mathematical definition of $A_O$ and, subsequently, $P_O$.
At first glance, our definitions may seem unnecessarily complex.
The benefit of such a rigorous approach is that it highlights the dependence of these quantities on the objectives $\mathbf{O}$.
Explicitly denoting this dependence will help us appreciate the nature of our approximation to the Pareto frontier when $\mathbf{O}$ is replaced by a Gaussian process surrogate.

Select a point $\mathbf{r}=(r_1,\dots,r_m)\in\mathbb{R}^m$ for which we have $\min_{\mathbf{x}\in X}O_i(\mathbf{x})\ge r_i$.
Since $X$ is compact, such a point exists if $O_i(\mathbf{x})$ is continuous.
$\mathbf{r}$ is known as the \emph{reference point}.
Consider the vector valued function $\mathbf{O}: X\rightarrow \mathbb{R}^m$ defined by $\mathbf{O} = (O_1,\dots,O_m)$.
$\mathbf{O}$ just joins all the expected objectives in a vector.
The image $\mathbf{O}[X]$ of $X$ under $\mathbf{O}$, defined by
$$
\mathbf{O}[X] = \left\{\by\in\mathbb{R}^m: \exists \mathbf{x}\in X, \by = \mathbf{O}(\bx) \right\},
$$
is the set of all achievable objectives.
We do not know exactly how $\mathbf{O}[X]$ looks like.
However, exploiting the definition of the reference point, we see that $\mathbf{O}[X]$ is fully contained in the $m$-dimensional cone $[\mathbf{r}, \infty] := \times_{i=1}^m[r_i, \infty)$, i.e.,
$$
    \mathbf{O}[X] \subset [\mathbf{r}, \infty).
$$
Consider any subset $B$ of $[\mathbf{r}, \infty)$.
    We define the attained set of $B$, denoted by $A[B]$, to be the set of points in $[\mathbf{r}, \infty)$ that are dominated by $B$, i.e.,
\begin{equation}
  \label{eqn:A}
  A[B] := \left\{\mathbf{y}\in [\mathbf{r},\infty): \exists \mathbf{y}'\in B, \mathbf{y}'\ge \mathbf{y}\right\},
\end{equation}
where $\by' \ge \by$ corresponds to element-wise comparison.
The attained set of our multi-objective problem is just:
\begin{equation}
    \label{eqn:A_moo}
    A_O := A[\mathbf{O}[X]].
\end{equation}
Finally, we define the Pareto frontier of $B$, denoted by $P[B]$, to be the set of points in $B$ that are not dominated by any other point in $B$, i.e.,
\begin{equation}
    P[B] := \left\{\by\in B: \{\by'\in B: \by' \ge \by\}=\emptyset\right\}.
\end{equation}
But we can get the Pareto frontier of $B$ directly from the boundary of its attained set.
Specifically, it is easy to prove that $P[B]$ is the top right boundary of $A[B]$, i.e.,
\begin{equation}
    \label{eqn:P}
    P[B] = \partial A[B] \setminus \cup_{i=1}^m\{\mathbf{r} + t(\max_{\by\in B}y_i)\mathbf{e}_i\},
\end{equation}
where $\mathbf{e}_i$ is the standard basis function of $\mathbb{R}^m$ pertaining to the $i$-th dimension.
The Pareto front of our multi-objective problem is just:
\begin{equation}
    \label{eqn:P_moo}
    P_O := P[\mathbf{O}[X]].
\end{equation}

Assume that we can choose to measure the QoIs at any design point $\bx\in X$ we wish, albeit only a limited number of times $n$.
Such measurements take place as follows.
When we request information about $\bx$, a latent process samples an \emph{unobserved} $\omega\in\Omega$ according to the probability measure $\mathbb{P}$, and we observe a noisy version of the QoIs $\by = (o_1(\bx,\omega),\dots,o_m(\bx,\omega))$.
This setup is general enough to account for both simulation-based and experiment-based QoIs.

Assume that we have queried the information source at $n$ design points.
\begin{equation}
  \bx_{1:n} = (\bx_1,\dots,\bx_n)\in X^n,
  \label{eqn:inputs}
\end{equation}
and that we have made the following noisy observations:
\begin{equation}
  \by_{1:n} = (\by_1,\dots,\by_n).
  \label{eqn:observed_outputs}
\end{equation}
We address two problems:
\begin{enumerate}
  \item What is our \emph{state of knowledge} about the true Pareto-efficient frontier $P_O$ given the observations $(\bx_{1:n},\by_{1:n})$?
  \item How should we select $\bx_{1:n}$ so that we come as close as possible to discovering the true Pareto-efficient frontier $P_O$?
\end{enumerate}
In the language of probability theory \fcite{jaynes2003}, the former problem seeks to characterize the probability (a state of belief) of a design being optimal conditional on the observations.
The uncertainty encoded in this probability is epistemic and it is induced by the fact that inference is based on just a small number of observations.
We address this problem by leveraging the Bayesian nature of Gaussian process surrogates, see \sref{gpr}.
Looking for an optimal information acquisition policy that solves the latter problem is a mathematically intractable task since the problem is equivalent to a non-linear stochastic dynamic program~\fcite{powell2012,bertsekas2007}.
We rely on a myopic/greedy one-step-look-ahead strategy (which is sub-optimal) by extending the definition of the standard EIHV, see \sref{EEIHV}, so that it can cope robustly with noise.

\subsection{Gaussian process regression}
\label{sec:gpr}

Gaussian process (GP) regression \fcite{rasmussen2006} is the Bayesian interpretation of classical Kriging \fcite{cressie1990,Smith2014}.
It is a powerful non-linear and non-parametric regression technique that is able to quantify the epistemic uncertainty induced by limited data.
We use GP regression to model our state of knowledge about the objectives, i.e., $O_i(\bx) = \E[o_i(\bx,\omega)],i=1,\dots,m$, as induced by a set of \emph{n} observations $(\bx_{1:n},\by_{1:n})$.
The methodology applies to each $i=1,\dots,m$ independently.
For simplicity, we will write $f(\bx)$ for $O_i(\bx)$ and $y_{1:n}$ for $y_{i,1:n}=(y_{i1},\dots,y_{in})$.

\subsubsection{Expressing prior beliefs}
\label{sec:prior}
Let $(\Omega^e,\mathcal{F}^e, \mathbb{P}^e)$ be the probability space corresponding to our epistemic uncertainty.
Note that this is different from $(\Omega,\mathcal{F},\mathbb{P})$ which is associated with the problem uncertainty.
A GP $f^e(x,\omega^e)$ is a $(\Omega^e,\mathcal{F}^e, \mathbb{P}^e)$-random field indexed by $\bx\in X$ with Gaussian finite dimensional distributions.
That is, for any $\bx_{1:n}\in X^n$ the random vector $f^e_{1:n}:=(f^e(\bx_1,\omega^e),\dots,f^e(\bx_n,\omega^e))$ follows a multivariate Gaussian.
The interpretation is as follows.
Nature has chosen a reality $\omega^e\in\Omega^e$, i.e., $f(\cdot)\equiv f^e(\cdot,\omega^e)$, that we cannot directly observe.
$(\Omega^e,\mathcal{F}^e,\mathbb{P}^e)$ models our prior state of knowledge about this reality, in the sense that for all $B\in\mathcal{F}^e$ the probability that we give to $\omega^e\in B$ is $\mathbb{P}^e[B]=\int_B d\mathbb{P}^e(\omega^e)$.

A GP is characterized by a mean and a covariance function.
Without loss of generality, we may assume that the mean function is zero, since the covariance can always be modified to include a non-zero mean trend.
Mathematically, we write:
\begin{equation}
  f^e |\theta^e \sim \GP(0, k),
\end{equation}
where $k:X\times X\times \Theta^e \rightarrow\R$ is a covariance function parameterized by the epistemic random variable $\theta^e:\Omega^e\rightarrow\Theta^e$.
According to the definition of the GP, our a priori beliefs about the values $f^e_{1:n}$ are captured by:
\begin{equation}
  \label{eqn:f_prior}
  f^e_{1:n}|\bx_{1:n},\theta^e \sim \calN(0, k(\bx_{1:n},\theta^e)),
\end{equation}
where $\mathcal{N}(\lambda,\Sigma)$ denotes the multivariate Gaussian distribution with mean $\lambda$ and covariance matrix $\Sigma$, for all $\bx'_{1:n'}\in X^{n'}$ we define $k(\bx_{1:n},\bx'_{1:n'}, \theta^e)$ to be the $n\times n'$ matrix with $(i,j)$ element $k(\bx_i, \bx_j,\theta^e)$, and $k(\bx_{1:n},\theta^e) := k(\bx_{1:n}, \bx_{1:n},\theta^e)$ is the covariance matrix.
In our numerical examples, we use the Matern($\nu=\frac{3}{2}$) \fcite{rasmussen2006} covariance:
\begin{equation}
        k(\bx,\bx',\theta^e) = {s^2} \Bigg( \exp\left\{ { - \sqrt{3\sum\limits_{j = 1}^d
  {\frac{{{{({x_j} - {x_j}')}^2}}}{{\ell_j^2}}}} }\right\} \Bigg) \Bigg( 1 + \sqrt{3\sum\limits_{j = 1}^d
  {\frac{{{{({x_j} - {x_j}')}^2}}}{{\ell_j^2}}}}\Bigg) ,
\end{equation}
where $d$ is the dimensionality of the design space, $s>0$ and $\ell_j>0$ can be interpreted as the signal strength of the response and the lengthscale along input dimension $j$, respectively, and $\theta^e  = (s,{\ell_1}, \ldots ,{\ell_d})\in\R^d_+$.
\subsubsection{Modeling the measurement process}
\label{sec:like}
In general, the noise that contaminates the measurement $y$ is heteroscedastic, i.e., input-dependent.
However, we approximate this noise as Gaussian with a fixed, but unknown, variance $\nu^2$.
Despite this fact, we observe numerically that the GP can still estimate the optimization objectives, i.e., expectation of $y$, when the noise to signal ratio is not too big.
The \emph{likelihood} of the model is:
\begin{equation}
 p(y_{1:n} | \bx_{1:n}, \theta^e) = \calN\left(y_{1:n} | 0, k(\bx_{1:n},\theta^e) + \nu^{2}I_n\right),
    \label{eqn:y_like}
\end{equation}
where $I_n\in\R^{n\times n}$ is the identity matrix,
$k(\bx_{1:n},\theta^e)$ is as in Eq.~(\ref{eqn:f_prior}),
and, for notational convenience, we have re-defined $\theta^e \leftarrow (\theta^e, \nu)$.

\subsubsection{Posterior state of knowledge about the objectives}
\label{sec:posterior}
Bayes rule combines our prior beliefs with the data
and yields a posterior probability measure on the space of meta-models.
Conditioned on the hyperparameters $\theta^e$, this measure is also a GP,
\begin{equation}
    \label{eqn:posterior_gp}
    f^e | \bx_{1:n}, y_{1:n}, \theta^e \sim \GP\left(\mu_n, k_n\right),
\end{equation}
where the posterior mean and covariance functions are
\begin{equation}
    \mu_n(\bx;\theta^e) = k_n(\bx,\bx_{1:n},\theta^e)\left[k(\bx_{1:n},\theta^e) + \nu^{2}I_n\right]^{-1}y_{1:n},
  \label{eqn:posterior_mean}
\end{equation}
and
\begin{equation}
  \begin{array}{ccc}
  k_n(\bx,\bx',\theta^e)&=&k(\bx, \bx',\theta^e) \\
             &&- k_n(\bx,\bx_{1:n},\theta^e)\left[k(\bx_{1:n},\theta^e) + \sigma^{2}I_n\right]^{-1}k_n(\bx_{1:n},\bx,\theta^e)\\
  \end{array}
\label{eqn:predictive_covariance}
\end{equation}
respectively.
Restricting our attention to a specific design point $\bx$,
we can derive from \qref{posterior_gp} the \emph{point-predictive} PDF conditioned on the hyperparameters $\theta^e$:
\begin{equation}
    \label{eqn:point_predictive}
    f^e(\bx) | \bx_{1:n}, y_{1:n}, \theta^e\sim
     \calN\left(\mu_n(\bx;\theta^e), \sigma_n^2(\bx;\theta^e)\right),
\end{equation}
where \emph{predictive variance} is $\sigma_n^2(\bx;\theta^e) = k_n(\bx,\bx;\theta^e)$.

The hyper-parameters of the covariance function are estimated by
maximizing the likelihood  $p(y_{1:n}|\bx_{1:n},\theta^e)$ with respect to $\theta^e$.
To avoid numerical instabilities, one typically works with the logarithm of the likelihood:
\begin{gather}
\mathcal{L}(\theta^e) =  - \frac{1}{2}{{y}_{1:n}}^T\left[k(\bx_{1:n},\theta^e) + \nu^{2}I_n\right]^{-1}{{y}_{1:n}} \nonumber\\
  - \frac{1}{2}\log \det\left[k(\bx_{1:n},\theta^e) + \nu^{2}I_n\right]- \frac{n}{2}\log 2\pi.
\label{eqn:mle}
\end{gather}
This maximization problem is solved using the BFGS algorithm
\fcite{Byrd1995}.
To account for the positivity constraints we simply optimize with respect to the logarithms of the hyperparameters.
The solution of this optimization problem, denoted by $\hat{\theta}^e$, is known as the maximum likelihood estimate (MLE) of $\theta^e$.
For notational convenience, in what follows we are not going to be explicitly indicating the dependence of $\mu_n$ and $k_n$ on $\theta^e$.
Instead it will be understood that $\mu_n(\bx) \equiv \mu_n(\bx,\hat{\theta}^e)$, $k_n(\bx,\bx')\equiv k_n(\bx,\bx',\hat{\theta}^e)$, and $\sigma_n(\bx)\equiv \sigma_n(\bx,\hat{\theta}^e)$.

\subsection{Characterization of the Pareto-efficient frontier using limited data}
\label{sec:pareto}
What is our state of knowledge about the true Pareto-efficient frontier $ P_O$ given $n\le N$ observations $(\bx_{1:n}, \by_{1:n})$?
Let $\mathbf{f}^e = (f^e_1,\dots, f^e_m)$ be the GPs representing our state of knowledge about each one of the $m$ objectives.
Our state of knowledge about the relation `$\succcurlyeq$' is now captured by the \emph{random relation} `$\succcurlyeq^e$', namely $\mathbf{x}\succcurlyeq^e\mathbf{x}'$ if and only if $\mathbf{f}^e(\bx) \ge \mathbf{f}^e(\mathbf{x}')$.
Our state of knowledge about the attained set $A_O$ of Eq.~(\ref{eqn:A_moo}) is given by the random set $A[\mathbf{f}^e[X]]$.
Similarly, our state of knowledge about the Pareto front $P_O$ of Eq.~(\ref{eqn:P_moo}) is represented by the random set $P[\mathbf{f}^e[X]]$.

The first step is to derive summary statistics of $A[\mathbf{f}^e[X]]$ that can be used to visualize our epistemic uncertainty about it.
Following \fcite{Binois:2015hg,chevalier2013}, we achieve this by estimating the Vorob'ev expectation and deviation of the random set $A[\mathbf{f}^e[X]]$.
Towards this end, we introduce the \emph{attainment function} and its upper level sets.
The attainment function $a^e_n:[\mathbf{r}, \infty)\rightarrow [0,1]$ is defined to be the conditional probability, given $(\mathbf{x}_{1:n},\mathbf{y}_{1:n})$, that a vector of objectives $\mathbf{y}\in [\mathbf{r},\infty)$ can be attained, i.e., we define
\begin{equation}
  \label{eqn:a_n}
  a^e_n(\mathbf{y}) := \mathbb{P}^e\left[\left\{\omega^e\in\Omega^e:\mathbf{y}\in A[\mathbf{f}^e_{\omega^e}[X]]\right\}|\mathbf{x}_{1:n},\mathbf{y}_{1:n}\right],
\end{equation}
where $\mathbf{f}^e_{\omega^e}(\cdot) = \left(f_1^e(\cdot,\omega^e),\dots,f_m^e(\cdot,\omega^e)\right)$.
For $\beta\in[0,1]$, the upper level sets of the attainment function,
\begin{equation}
  \label{eqn:Q_beta}
  Q_{n,\beta}^e := \left\{\mathbf{y}\in[\mathbf{r},\infty): a^e_n(\mathbf{y})\ge \beta\right\},
\end{equation}
are known as the $\beta$-quantiles of $ A[\mathbf{f}^e[X]]$. Intuitively, $Q_{n,\beta^*}^e$ can be seen as the set of objectives that are considered achievable with probability greater than or equal to $\beta$.
The conditional \emph{Vorob'ev expectation} \fcite{molchanov2005} of $A[\mathbf{f}^e[X]]$ is defined to be the $\beta^*$-quantile $Q_{n,\beta^*}^e$ for which:
\begin{equation}
  \label{eqn:Vorobev}
  \lambda(Q_{n,\beta}^e) \le \mathbb{E}^e\left[\lambda\left(A[\mathbf{f}^e[X]]\right)|\mathbf{x}_{1:n},\mathbf{y}_{1:n}\right]\le \lambda(Q_{n,\beta^*}^e),\;\forall\beta \in[\beta^*,1],
\end{equation}
where $\lambda$ is the Lebesgue measure on $\mathbb{R}^m$.
In words, $Q_{n,\beta^*}^e$ is the $\beta$-quantile that has the same Lebesgue measure as the conditional expectation of the Lebesgue measure of the attained set.
Intuitively, $Q_{n,\beta^*}^e$ and its top right boundary are our expectations about the attained set $A_O$ and $P_O$, respectively, after observing $(\mathbf{x}_{1:n},\mathbf{y}_{1:n})$.

Now, we are in a position to quantify our uncertainty about $P_O$.
Consider the symmetric difference $Q_{n,\beta^*}^e\triangle A[\mathbf{f}^e[X]]$ between the set $Q_{n,\beta^*}^e$ and $A[\mathbf{f}^e[X]]$ defined by
\begin{equation}
  \label{eqn:sym_diff}
  Q_{n,\beta^*}^e\triangle A[\mathbf{f}^e[X]] := \left(Q_{n,\beta^*}^e\cup A[\mathbf{f}^e[X]]\right)\setminus\left(Q_{n,\beta^*}^e\cap A[\mathbf{f}^e[X]]\right).
\end{equation}
That is, a point $\mathbf{y}$ belongs in $Q_{n,\beta^*}^e\triangle A[\mathbf{f}^e[X]]$ only if it belongs to exactly one of these sets.
Such points appear in the top right corner of $[\mathbf{r},\infty)$ and are candidate points for the Pareto front.
    Therefore, we quantify our uncertainty about $P_O$ through the \emph{symmetric deviation function} $d_{n}^e:[\mathbf{r},\infty)\rightarrow [0,1]$ defined as the conditional probability that a vector of objectives $\mathbf{y}\in[\mathbf{r},\infty)$ belongs to the symmetric difference $Q_{n,\beta^*}^e\triangle A[\mathbf{f}^e[X]]$, i.e.,
\begin{equation}
  \label{eqn:sym_dev}
  d_{n}^e(\mathbf{y}) = \mathbb{P}^e\left[\mathbf{y}\in Q_{n,\beta^*}^e\triangle A[\mathbf{f}^e[X]]|\mathbf{x}_{1:n},\mathbf{y}_{1:n}\right].
\end{equation}

Unfortunately, it is not possible to characterize $a_n^e(\mathbf{y})$, $Q_{n,\beta^*}^e$, and $d_n^e(\mathbf{y})$ exactly.
The difficulty arises from the fact that $X$ may be infinite dimensional.
To overcome this obstacle, we use a Monte Carlo (MC) approach.
Let $(\tilde{\Omega}, \tilde{\mathcal{F}}, \tilde{\mathbb{P}})$ be a new probability space associated with the MC approximation uncertainty.
Let $\tilde{X}_s:\tilde{\Omega}\rightarrow X^{\tilde{n}}$, collectively denoted by $\tilde{X}_{1:S} = (\tilde{X}_1,\dots,\tilde{X}_S)$, be independent identically distributed (iid) random variables in $(\tilde{\Omega}, \tilde{\mathcal{F}}, \tilde{\mathbb{P}})$ with values in $X^{\tilde{n}}$.
For each one, we have $\tilde{X}_s := \tilde{\mathbf{x}}_{1,1:\tilde{n}} := (\tilde{\mathbf{x}}_{s,1},\dots,\tilde{\mathbf{x}}_{s,\tilde{n}})$
The specific distribution of these variables is not important as soon they cover $X$.
For convergence, it suffices to make all the $\tilde{\mathbf{x}}_{s,i},s=1,\dots,S,i=1,\dots,\tilde{n}$ iid with a support that covers $X$.
In our numerical examples, we take all these random variables to be independently uniform.
Conditional on each $\tilde{X_s}$, define the epistemic random variable $\tilde{F}_s^e\in\mathbb{R}^{m\tilde{n}}$ associated with the values of the objectives on $\tilde{X}_s$.
That is, $\tilde{F}_s^e:=\tilde{\mathbf{f}}^e_{s,1:m,1:\tilde{n}}:= (\tilde{f}^e_{s,1,1:\tilde{n}},\dots,\tilde{f}^e_{s,m,1:\tilde{n}})$, with $\tilde{f}_{s,i,1:\tilde{n}} := (f^e_i(\tilde{\mathbf{x}}_{s,1}),\dots,f^e_i(\tilde{\mathbf{x}}_{s,\tilde{n}}))\in\mathbb{R}^{\tilde{n}}$.
Note that, since we constructed each one of the GPs representing the objectives independently, we have that $\tilde{f}^e_{s,i,1:\tilde{n}},s=1,\dots,S,i=1,\dots,m$ are independent.
Making use of the posterior GP representing our state of knowledge about $f^e_i(\mathbf{x})$, see Eq.~(\ref{eqn:posterior_gp}), we get that, conditional on $\tilde{\mathbf{x}}_{1:\tilde{n}}$ and $(\mathbf{x}_{1:n}, y_{i,1:n})$, $\tilde{f}^e_{s,i,1:\tilde{n}}$ is normally distributed:
\begin{equation}
  \tilde{f}^e_{s,i,1:\tilde{n}} | \tilde{\mathbf{x}}_{s,1:\tilde{n}},\mathbf{x}_{1:n},y_{i,1:n} \sim \mathcal{N}\left(\mu_{i,n}(\tilde{\mathbf{x}}_{s,1:\tilde{n}}), k_{i,n}(\tilde{\mathbf{x}}_{s,1:\tilde{n}})\right),
\end{equation}
where $\mu_{i,n}(\mathbf{x})$ and $k_{i,n}(\mathbf{x},\mathbf{x}')$ are the posterior mean and posterior covariance functions ($\mu_n(\mathbf{x})$ and $k_n(\mathbf{x},\mathbf{x}')$) of Sec.~\ref{sec:posterior}, respectively, if we make the substitution $y_{1:n} \leftarrow y_{i,1:n}$.
Using $\tilde{F}^e_s$, and the definition in \qref{A} we denote the sampled attained set by $A[\tilde{F}^e_s]$ and the
corresponding \emph{sampled Pareto front} by $P[\tilde{F}^e_s]$.
Now we can compute the \emph{empirical attainment function} $\tilde{a}_{S,\tilde{n},n}^e:[\mathbf{r},\infty)\rightarrow [0,1]$:
\begin{equation}
  \label{eqn:empirical_a}
  \tilde{a}_{S,\tilde{n},n}^e(\mathbf{y}) = \frac{1}{S}\sum_{s=1}^S 1_{A[\tilde{F}^e_s]}(\mathbf{y}),
\end{equation}
where $1_B(\mathbf{y})$ is the characteristic function of the set $B$.
Using $\tilde{a}_{S,\tilde{n},n}^e(\mathbf{y})$ we can obtain estimates of the $\beta$-quantiles, say $\tilde{Q}_{S,\tilde{n},n,\beta}^e$.
Just like \fcite{Binois:2015hg}, estimates of the $\beta$-quantiles can be used within a bisection algorithm to estimate the Vorob'ev expectation $\tilde{Q}_{S,\tilde{n},n,\beta^*}^e$.
Finally, we compute the \emph{empirical symmetric deviation function}:
\begin{equation}
  \label{eqn:empirical_d}
  \tilde{d}_{S,\tilde{n},n}^e(\mathbf{y}) = \frac{1}{S}\sum_{s=1}^S1_{\tilde{Q}_{S,\tilde{n},n,\beta^*}\triangle A[\tilde{F}^e_s]}(\mathbf{y}),
\end{equation}
which is an estimate of $d_n^e(\mathbf{y})$.
In our numerical examples (in which $m=2$) we represent $\tilde{a}^e_{S,\tilde{n},n}(\mathbf{y})$ and $\tilde{d}_{S,\tilde{n},n}^e(\mathbf{y})$ on a $64\times 64$ grid defined on $\times_{i=1}^m[r_i,u_i]$, where $\mathbf{u} = (u_1,\dots,u_m)\in\mathbb{R}^m$ is a point of the design space with $u_i \ge \max_{\mathbf{x}\in X} O_i(\mathbf{x}),i=1,\dots,m$.
For larger number of objectives $m > 3$, more sophisticated techniques must be developed in order to overcome the curse of dimensionality.
From the law of large numbers, we have that
\begin{equation}
  \lim_{S\rightarrow\infty}\lim_{\tilde{n}\rightarrow\infty} \tilde{a}^e_{S,\tilde{n},n} = a^e_n,
\end{equation}
\begin{equation}
  \lim_{S\rightarrow\infty}\lim_{\tilde{n}\rightarrow\infty}\tilde{d}^e_{S,\tilde{n},n} = d^e_n.
\end{equation}
We also expect that the attainment function $a^e_n$ will converge to the characteristic function of the attained set $A_O$ as $n\rightarrow \infty$ on a set of design points that becomes dense.
The exact nature of the latter convergence is beyond the scope of the present work.

\subsection{Extended expected improvement over dominated hypervolume}
\label{sec:EEIHV}
Given our current state of knowledge about $P_O$, how should we select the next observation $\bx$?
We derive a myopic one-step-look-ahead strategy that attempts to sequentially maximize the expected improvement in the volume of the attained set.
Specifically, we define the \emph{extended expected improvement over the dominated hypervolume} (EEIHV) as the expectation of the change in the Lebesgue measure of the attained set conditional on a hypothetical observation.
Mathematically, we define for $\bx\in X$:
\begin{equation}
  \label{eqn:eihv}
  \begin{array}{cccl}
      \operatorname{EEIHV}(\bx) &=& 
      \mathbb{E}^e\Big[&
  \mathbb{E}^e\big[\lambda\left(A[\mathbf{f}^e[X]]\right)\big|\mathbf{x},\mathbf{y},\mathbf{x}_{1:n},\mathbf{y}_{1:n}\big]\\
  &&&    -\mathbb{E}^e\big[\lambda\left(A[\mathbf{f}^e[X]\right)|\mathbf{x}_{1:n},\mathbf{y}_{1:n}\big]
\Big|\mathbf{x},\mathbf{x}_{1:n},\mathbf{y}_{1:n}\Big],
\end{array}
\end{equation}
where the outer expectation is over our state of knowledge about the hypothetical measurement $\by$ induced by the GPs of \sref{gpr}:
\begin{equation}
  \label{eqn:p_y}
  p(\by|\bx,\bx_{1:n},\by_{1:n}) = \prod_{i=1}^{m}\calN(y_{i}| \mu_{i,n}(\bx), \sigma_{i,n}^2(\bx;\theta^e) + \nu^2),
\end{equation}
where $\mu_{i,n}(\cdot) = \mu_{i,n}(\cdot;\theta^e_i)$ and $\sigma_{i,n}^2(\cdot)=\sigma_{i,n}^2(\cdot;\theta^e_i)$ are the posterior predictive mean and variance of the GP $f^e_i$ pertaining to objective $i=1,\dots,m$, see \qref{point_predictive}.
Our myopic strategy is outlined in Algorithm~\ref{alg:bgo}.

\qref{eihv} is analytically intractable and must be approximated using the sampling methods of \sref{pareto}. 
This is computationally inefficient because it does not allow the use of gradient-based optimization algorithms such as BFGS.
To overcome this difficulty, we derive an approximation that will allow us to make use of the analytical formulas derived by \fcite{emmerich2008computation}.
We have:
\begin{eqnarray*}
    \mathbb{E}^e\big[\lambda\left(A[\mathbf{f}^e[X]]\right)\big|\mathbf{x}_{1:n},\mathbf{y}_{1:n}\big] &\ge&
     \mathbb{E}^e\big[\lambda\left(A[\mathbf{f}^e[\bx_{1:n}]]\right)|\mathbf{x}_{1:n},\mathbf{y}_{1:n}\big]\\
     &\approx&
    \lambda\left(A[\boldsymbol{\mu}_n[\bx_{1:n}]]\right).
\end{eqnarray*}
The first row inequality comes from $\bx_{1:n}\subset X$ implying $\mathbf{f}^e[\bx_{1:n}]\subset \mathbf{f}^e[X]$ which, in turn, yields $A[\mathbf{f}^e[\bx_{1:n}]]\subset A[\mathbf{f}^e[X]]$. 
For the approximation in the second row, start by noticing that $\mathbf{z}=\mathbf{f}^e[\bx_{1:n}]$ conditioned on $\bx_{1:n}$ and that $\by_{1:n}$ follows a multivariate Gaussian, see \qref{posterior_gp}.
Then, take the Taylor expansion of $\lambda(A[\mathbf{z}])$ about $\mathbf{z}=\mathbf{z}_0=\boldsymbol{\mu}_n(\bx_{1:n}):=\left(\mu_{1,n}(\bx_{1:n}),\dots,\mu_{m,n}(\bx_{1:n})\right)$.
The zero order term is the constant you see above, i.e., $\lambda\left(A[\boldsymbol{\mu}_n[\bx_{1:n}]]\right)$.
The expectation of the first order term vanishes and we ignore second and higher order terms.
Thinking in the same way, we can get:
\begin{eqnarray*}
        \mathbb{E}^e\big[
            \lambda\left(A[\mathbf{f}^e[X]]\right)
            \big| \bx, \by, \bx_{1:n}, \by_{1:n}\big] &\ge&
            \mathbb{E}^e\big[
                \lambda\left(A[\mathbf{f}^e[\bx_{1:n}\cup \{\bx\}]]\right)
            \big| \bx, \by, \bx_{1:n}, \by_{1:n}\big]\\
        &\approx& \lambda\left(A[\boldsymbol{\mu}_{n,(\bx,\by)}[\bx_{1:n}\cup\{\bx\}]]\right),
\end{eqnarray*}
where $\boldsymbol{\mu}_{n,(\bx,\by)}$ is the posterior mean after seeing the hypothetical observation $(\bx,\by)$.
Finally, we approximate the expectation over the hypothetical measurement as:
\begin{eqnarray*}
    &&\mathbb{E}^e\Big[
        \lambda\left(A[\boldsymbol{\mu}_{n,(\bx,\by)}[\bx_{1:n}\cup\{\bx\}]]\right)
        \Big|\bx, \bx_{1:n}, \by_{1:n}\Big] \approx\\
    &&\mathbb{E}^e\Big[
    \lambda\left(A[\boldsymbol{\mu}_{n,\left(\bx,\boldsymbol{f}^e(\bx)\right)}[\bx_{1:n}\cup\{\bx\}]]\right)
\Big|\bx, \bx_{1:n}, \by_{1:n}\Big].
\end{eqnarray*}
To see why this is possible, note that $\by = \mathbf{f}^e(\bx) + \nu^2 \boldsymbol{\epsilon}$ where $\boldsymbol{\epsilon}$ is Gaussian with zero mean and unit covariance, take the Taylor expansion of the integrand in the first line about $\boldsymbol{\epsilon}=\mathbf{0}$, and keep only the zero order term (the expectation of the first order term vanishes).
Putting everything together, we get the (approximate) inequality:
\begin{equation}
    \label{eqn:EEIHV}
    \begin{array}{ccc}
    \operatorname{EEIHV}(\bx) &\tilde{\ge}& \overline{\operatorname{EEIHV}}(\bx) \\ 
        &:=& \mathbb{E}^e\Big[
    \lambda\left(A[\boldsymbol{\mu}_{n,\left(\bx,\boldsymbol{f}^e(\bx)\right)}[\bx_{1:n}\cup\{\bx\}]]\right)
\Big|\bx, \bx_{1:n}, \by_{1:n}\Big]\\
    && - \lambda\left(A[\boldsymbol{\mu}_n[\bx_{1:n}]]\right).
    \end{array}
\end{equation}
The inequality is approximate because the first term on the right hand side is approximately greater than the second one.
The accuracy is again second order and proving it requires taking the Taylor expansion of the integrand of the first term with respect to $\mathbf{z}\equiv\mathbf{f}^e(\bx)$ about $\mathbf{z}=\mathbf{z}_0\equiv\boldsymbol{\mu}_n(\bx)$.

The important observation here is that the lower bound to EEIHV, i.e., $\overline{\operatorname{EEIHV}}$ on right hand side of \qref{EEIHV}, is similar to the original EIHV of \fcite{emmerich2008computation} with a few key differences.
Specifically, $\overline{\operatorname{EEIHV}}$ has the same analytical form as EIHV if in EIHV (i) we replace the observed targets with their projections to the posterior mean, i.e., if we work with the denoised measurements instead of the noisy ones; and (ii) we remove the noise variance from the predictive distribution of the GP.
Therefore, the analytical formula for the calculation of EIHV found in \fcite{emmerich2008computation} applies to $\overline{\operatorname{EEIHV}}$ subject to the aforementioned substitutions.
In all our numerical examples, we use $\overline{\operatorname{EEIHV}}$.
We maximize the lower bound over $\bx$ using BFGS with multiple random restarts.

\begin{algorithm}[h!]
    \caption{Information acquisition strategy for discovering the Pareto-frontier.}
\begin{algorithmic}[1]
\Require Initially observed designs ${\bf{x}}_{1:n}$;
         Initial objective measurements $\by_{1:n}$;
         number of restarts of EEIHV optimization $n_d$;
         maximum number of allowed information source queries $N_{\max}$;
         EEIHV tolerance  $\delta>0$.
            \While {$n < N_{\max}$}
                \State Train the GP for each objective as described in \sref{gpr}.
                \State Find $\bx_{n+1} = \arg\max_{\bx \in X}\overline{\operatorname{EEIHV}}(\bx)$ using $n_d$ random restarts of BFGS.
                \If{$\overline{\operatorname{EEIHV}}(\bx_{n+1}) < \delta$}
                    \State Break.
                \EndIf
                \State Evaluate the objectives at $ \bx_{n+1}$ measuring $ \by_{n+1}$.
                \State $\bx_{1:n+1} \leftarrow (\bx_{1:n},\bx_{n+1})$.
                \State $\by_{1:n+1} \leftarrow (\by_{1:n},\by_{n+1})$.
                \State $n\leftarrow n + 1$.
  \EndWhile
\end{algorithmic}
\label{alg:bgo}
\end{algorithm}

\input{results}

\input{wire-prob}

\input{conclusions}

\section{Acknowledgments}
Ilias Bilionis acknowledges the startup support provided by the School of 
Mechanical Engineering at Purdue University.\\
The authors acknowledge the support provided by Tata Consultancy Services, 
Pune, India.
\nolinenumbers
\bibliography{references}
\bibliographystyle{abbrv}
\end{document}

%% file: abs.tex
\begin{abstract}
\textit{
Design optimization of engineering systems with multiple competing objectives is a painstakingly tedious process especially when the objective functions are expensive-to-evaluate computer codes with parametric uncertainties.
The effectiveness of the state-of-the-art techniques is greatly diminished because they require a large number of objective evaluations, which makes them impractical for problems of the above kind. Bayesian global optimization (BGO), has managed to deal with these challenges in solving single-objective optimization problems and has recently been extended to multi-objective optimization (MOO).
BGO models the objectives via probabilistic surrogates and uses the epistemic uncertainty to define an information acquisition function (IAF) that quantifies the merit of evaluating the objective at new designs.
This iterative data acquisition process continues until a stopping criterion is met.
The most commonly used IAF for MOO is the expected improvement over the dominated hypervolume (EIHV) which in its original form is unable to deal with parametric uncertainties or  measurement noise.
In this work, we provide a systematic reformulation of EIHV to deal with 
stochastic MOO problems.
The primary contribution of this paper lies in being able to filter out the noise and reformulate the EIHV without having to observe or estimate the stochastic parameters.
An addendum of the probabilistic nature of our methodology is that it enables us to characterize our confidence about the predicted Pareto front.
We verify and validate the proposed methodology by applying it to synthetic test problems with known solutions.
We demonstrate our approach on an industrial problem of die pass design for a steel wire drawing process.
}
\end{abstract}

%% file: intro.tex
\section{Introduction}
\label{sec:intro}
The goal of this paper is to derive a sequential information acquisition 
methodology that aims at efficiently discovering the Pareto set of a 
stochastic MOO problem.
Stochastic MOOs are characterized by uncertain objective measurements, i.e.,
for a fixed design, repeated measurements of the objectives may vary.
When the objectives are the outcomes of an experiment, this randomness may be due
to manufacturing imperfections, operational uncertainties, wear and tear of the
specimen, sensor malfunction, etc.
When the objectives depend on a simulation model, then this randomness may be
induced by uncertainty in the model parameters, e.g., boundary/initial conditions,
parameters of constitutive relations, or artifact geometries.
In the latter case, the designer chooses probability distributions for all
uncertain parameters in an effort to accurately describe their state of
knowledge about the artifact.

MOO techniques based on evolutionary algorithms \fcite{deb2008}, e.g., the strength Pareto evolutionary algorithm \fcite{Zitzler2001}, the non-dominated sorting genetic algorithm II (NSGA-II) \fcite{Deb2002}, require a significant number of objective evaluations, especially when coupled with a sample average approximation \fcite{shapiro2003} to estimate the stochastic objectives. Other popular techniques like \emph{goal programming} \fcite{charnes1977,marler2004} that involve a slight modification of the original MOO objectives face shortcomings \fcite{zeleny1981} like selecting the relative importance of the objectives, or requiring the designer to have prior information about discontinuities in the objective space.

Bayesian global optimization (BGO) \fcite{mockus2012bayesian,jones1998} is a class of black-box optimization algorithms that can operate under a limited objective evaluation budget.
BGO models the objectives using probabilistic surrogates, e.g., Gaussian process regression, and exploits the epistemic uncertainty to select which experiments/simulations to perform.
The latter is typically done by maximizing an information acquisition function (IAF) which quantifies the value of evaluating the objective at a specific design.
The choice of the IAF depends on the details of the underlying optimization task.
One of the most popular IAFs is the expected improvement (EI) \fcite{mockus1975,jones1998,huang2006,ginsbourger2007multi,brochu2010}. 
The EI balances the exploration-exploitation trade-off better than other popular IAFs such as the probability of improvement (PI) or the upper confidence bound (UCB) \fcite{jones2001}.
Keane \fcite{keane2006} extended the original version of EI to MOO by deriving the expected improvement over the dominated hypervolume (EIHV).
The EIHV evaluates the expected improvement in the volume of the attained set induced by a hypothetical observation at an untried design.
\cite{emmerich2008computation} derived a closed form representation which made the evaluation of EIHV computationally efficient.
Research in EIHV has been gaining momentum over the past few years~\fcite{tesch2013,shimoyama2013,feliot2015}, but it has not yet been extended to cover the case of stochastic multi-objective optimization.

In this work, we propose an extension to the EIHV suitable for stochastic MOO, which is
the main contribution of this paper. 
We will be referring to the proposed methodology as the \emph{extended EIHV} (EEIHV).
Our proposal is a generalization of the extended expected improvement (EEI) which we 
developed in \fcite{Pandita2016} to deal with stochastic single-objective optimization.
The methodology relies on building probabilistic surrogates of the objectives 
and uses the EEIHV IAF to quantify the merit of evaluating the expensive stochastic computer 
code at a new design.
We leverage the work done in \fcite{Binois:2015hg} to 
quantify our uncertainty about the estimated PF at each stage/iteration.

We apply the above methodology to solve a multi-pass steel wire manufacturing
problem under uncertainty.
The competing objectives in this problem are the ultimate tensile strength (UTS) 
and the strain non-uniformity factor (SNUF) of the drawn wire. A  finite element (FE) solver
(developed at Tata Consulatancy Services (TCS), Pune, India) generates these objectives.
The reduction ratios and the die angles at each pass are the design/process variables 
which have associated uncertainties due to unavoidable manufacturing tolerances 
as well as die wear during the process.

The outline of the paper is as follows.
We start \sref{metho} by providing the mathematical definition of the stochastic 
MOO optimization problem that we are studying.
In \sref{gpr}, we introduce Gaussian process regression (GPR) which is used to 
construct the probabilistic surrogates of the map between the design variables and the objectives.
In \sref{EEIHV}, we derive our extension to EIHV suitable for stochastic multi-objective optimization.
Our numerical results are presented in \sref{results}.
In particular, in Sec.~\ref{sec:validate_2d} and~\ref{sec:validate_6d}, we validate 
our approach using two synthetic stochastic MOO problems with known analytical expressions, 
and we experiment with varying levels of stochasticity (to represent noisy measurements).
In \sref{wmp}, we apply our methodology to solve the wire drawing problem.
We present our conclusions in \sref{conclusions}.

%% file: results.tex
\section{Numerical Results}
\label{sec:results}
In \sref{visualize} we use a synthetic example to visualize some of the concepts used through out this section.
In Sections~\ref{sec:validate_2d} and~\ref{sec:validate_6d}, we validate our approach using two synthetic stochastic optimization problems with known optimal
solutions.
To assess the robustness of the methodology, we experiment with various levels
of stochasticity which causes the resultant
noise in the outputs. In \sref{wmp}, we solve the steel wire drawing 
problem with uncertainties in the incoming wire diameters and the die angles at 
each pass.
In all the problems the objectives are scaled by subtracting and dividing by the emprical mean and standard deviation, respectively.
\begin{figure}[h!]
    \centering
    \subfigure[]{
        \includegraphics[width=0.5\textwidth]{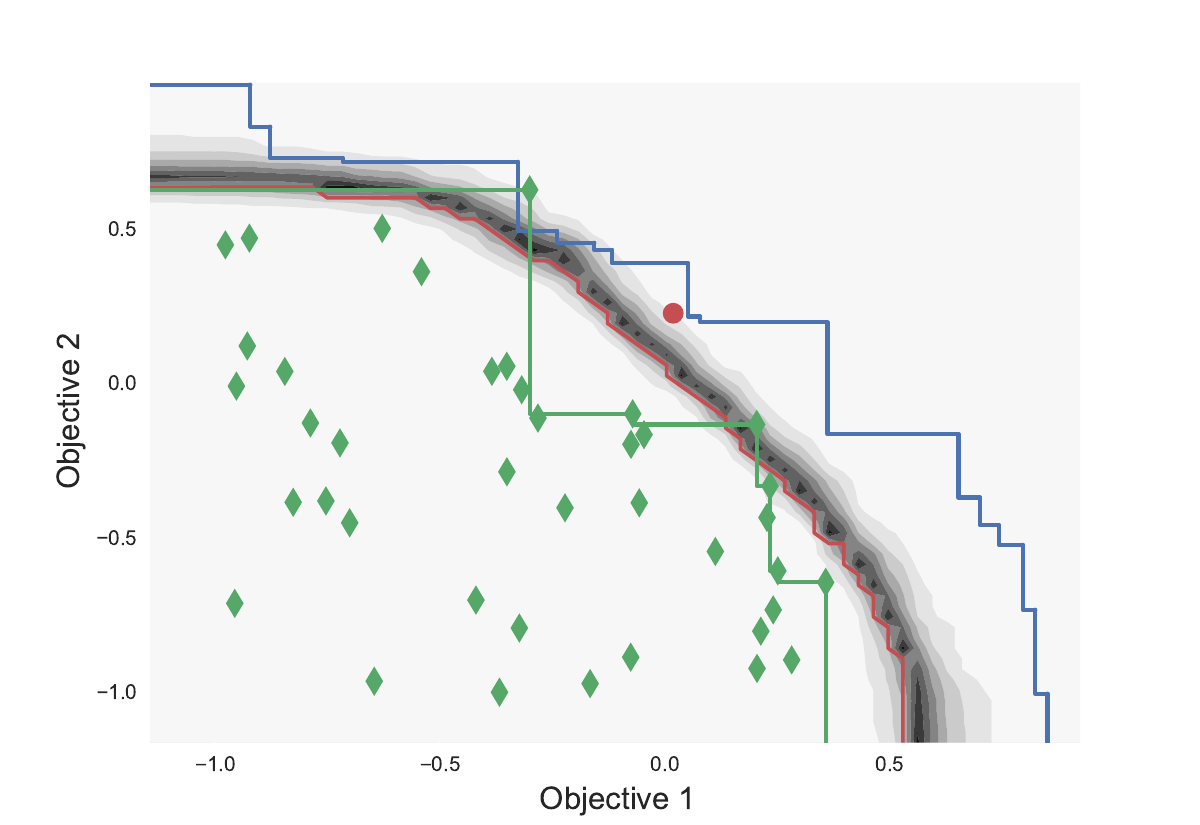}
    }
    \caption{A synthetic example of the template followed throughout
        the paper depicting the Pareto front and the
        representation of the uncertainty around the Pareto front. 
    }
    \label{fig:state-1}
\end{figure}

\subsection{Correspondence between nomenclature and visualizations}
\label{sec:visualize}
\fref{state-1} uses an $m=2$ synthetic example to help us visualize and name some of the concepts used throughout this section.
The dark blue staircase is an approximation of the true $P_O$, generated by taking the empirical Pareto frontier of sample averaged objective measurements at a large number of designs.
The figure also shows a scatter plot of the denoised measurements $\boldsymbol{\mu}_n(\bx_{1:n})$ (green dots), and as well as the corresponding empirical Pareto frontier $P[\boldsymbol{\mu}_n[\bx_{1:n}]]$ (green line).
The red dot marks the denoised measurement made at the design $\bx_{n+1}$ that maximizes $\overline{\operatorname{EEIHV}}(\bx)$.
The red line is the top right boundary of the Vorob'ev expectation of $A[\mathbf{f}^e[X]]$ conditioned on the observed data $(\bx_{1:n}, \by_{1:n})$, i.e., it is our expectation about $P[\mathbf{f}^e[X]]$ conditioned on our current state of knowledge.
The gray contours show to the symmetric deviation $d^e_n(\by)$ of $A[\mathbf{f}^e[X]]$ which corresponds to our uncertainty about $P[\mathbf{f}^e[X]]$.

\subsection{Two-dimensional synthetic example}
\label{sec:validate_2d}
Consider the two-dimensional synthetic multi-objective problem taken from \cite{parr2013}
which has been slightly modified for our use here:
\begin{gather}
    o_{1}(\bx,\omega) = -\bigg({b_{2} - {\frac{5.1}{4{\pi^{2}}}}{b_{1}}^{2} + {{\frac{5}{\pi}}b_{1} - 6}}\bigg)^{2} - \\
                      10\bigg[\bigg(1-\frac{1}{8\pi}\bigg)\cos(b_{1}) + 1\bigg], \nonumber\\
    o_{2}(\bx,\omega) = \sqrt{|(10.5-b_{1})||(b_{1} + 5.5)||(b_{2} + 0.5)|}+ \\
                      \frac{1}{30}\bigg(b_{2} - {\frac{5.1}{4{\pi^{2}}}}{b_{1}}^{2} - 6 \bigg)^{2}+ \nonumber\\
                      \frac{1}{3}\bigg[\bigg(1-\frac{1}{8\pi}\bigg)\cos(b_{1}) + 1\bigg], \nonumber\\
    b_{1}(\bx,\omega) = 15(x_{1} +s\xi(\omega)) - 5,\\
    b_{2}(\bx,\omega) = 15(x_{2} + s\xi(\omega)),
    \label{eqn:2d_example}
\end{gather}
for $\bx = (x_1,x_2) \in X = [0,1]^2$.
The $(\Omega,\mathbb{P},\mathcal{F})$ random variable $\xi$ is standard normal, i.e., $\xi\sim \mathcal{N}(0,1)$.
The parameter $s$ controls the standard deviation of the noise infused by $\xi$.
Notice that even though $\xi$ is normal, the measured objectives $o_i(\bx,\omega)$ are not normally distributed due to the non-linearities.
That is the statistics of the measurement process do not match our assumptions in \sref{gpr}.
We do this on purpose.
In real applications the statistics of the measurements process are not known and we would like to investigate to what extent the normality assumption produces robust results.

To validate our methodology, we must first estimate accurately the true $P_O$.
We achieve this by finding the empirical Pareto frontier of a large number of designs (10000)  while approximating $O_i(\bx)=\E[o_i(\bx,\omega)]$ with 100 Monte Carlo samples.
In this example, we aim to maximize the two objectives.

\begin{figure}[h!]
    \subfigure[]{
        \includegraphics[width=0.5\textwidth]{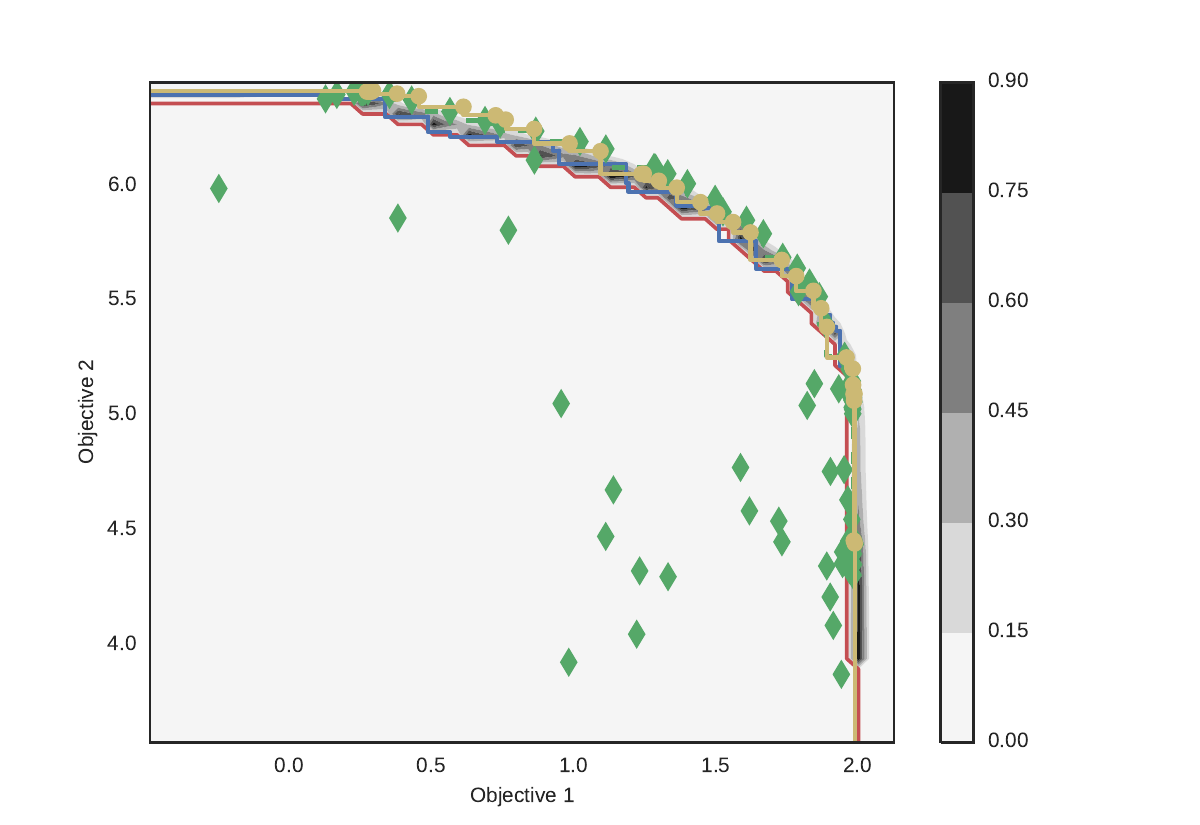}
    }
    \subfigure[]{
        \includegraphics[width=0.5\textwidth]{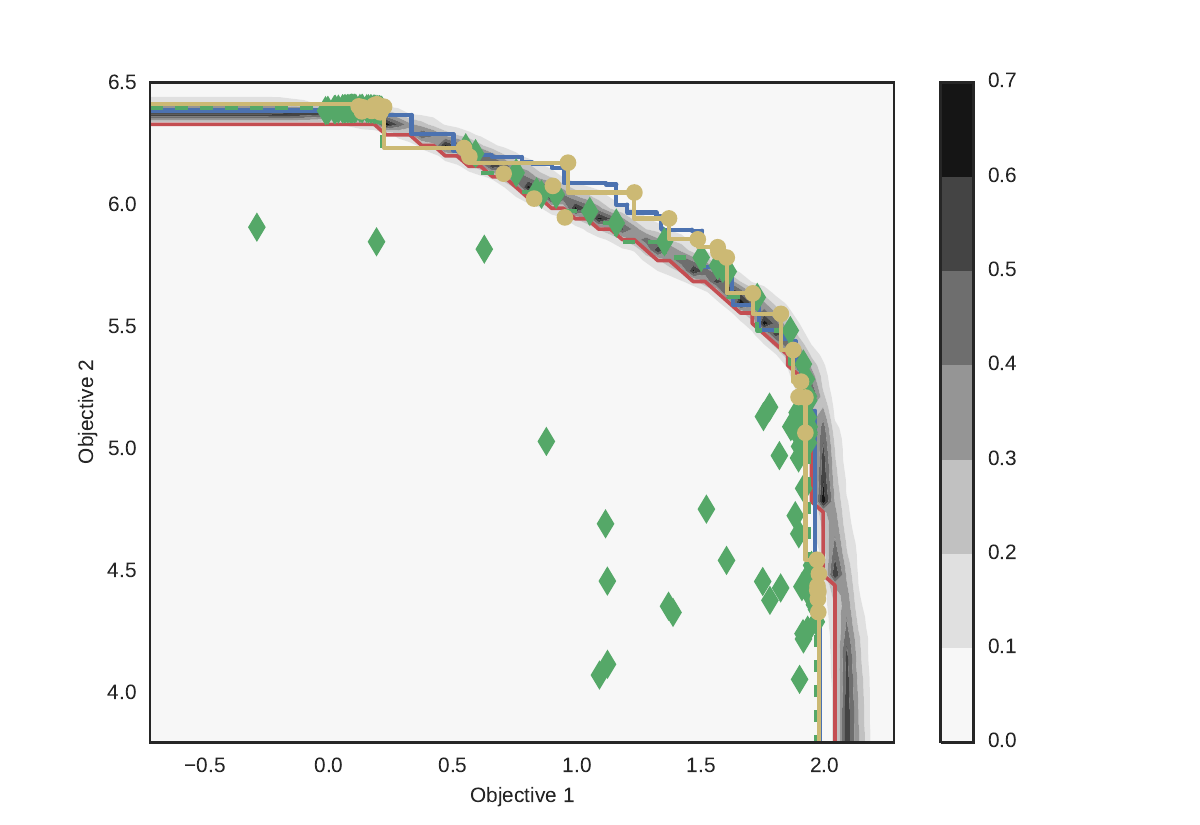}
    }
    \subfigure[]{
        \includegraphics[width=0.5\textwidth]{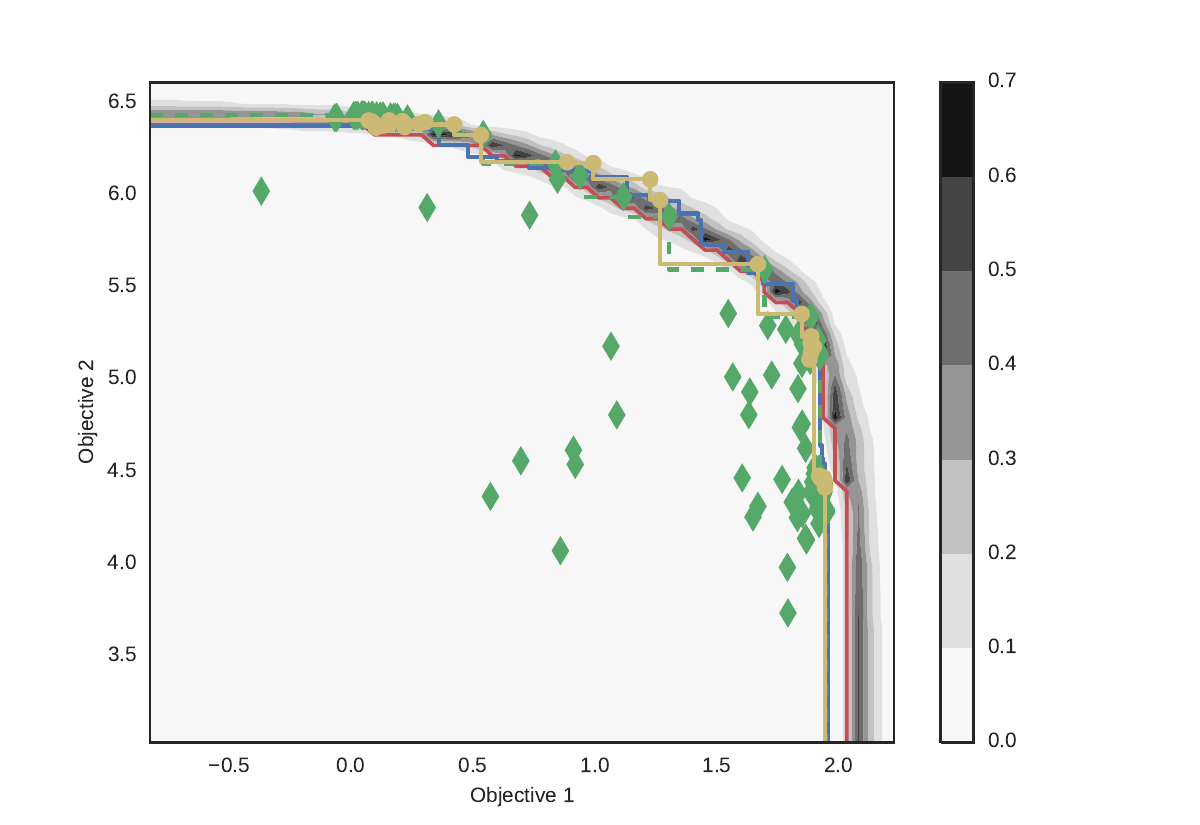}
    }
    \subfigure[]{
        \includegraphics[width=0.5\textwidth]{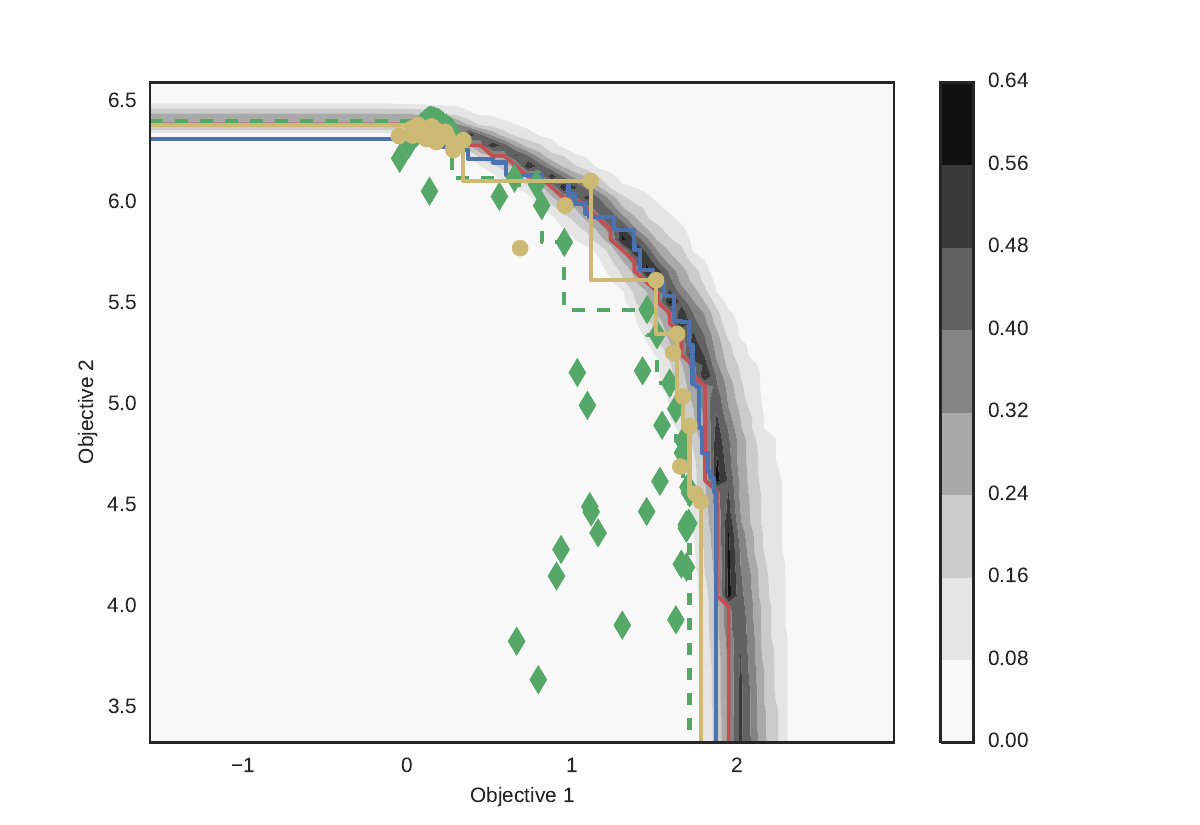}
    }
    \caption{Two-dimensional synthetic example for starting from $n=20$ initial measurements.
        Subfigures~(a)~($s=0.01$), ~(b)~($s=0.03$), ~(c)~($s=0.05$), and~(d)~($s=0.1$), 
        depict our state of knowledge about the final
        $P[\mathbf{f}^e[\bx_{1:n}]]$ after 100 measurements selected using Algorithm~\ref{alg:bgo}.
        }
    \label{fig:ex1_a}
\end{figure}

We start with $n=20$ random initial observations and we add an additional 100 measurements selected according to Algorithm~\ref{alg:bgo}.
\fref{ex1_a} depicts our final state of knowledge about $P[\mathbf{f}^e[\bx_{1:n}]]$ for increasing noise levels $s = 0.01, 0.03, 0.05,$ and $0.1$.
Another graphic that appears on this figure is the line joining the large yellow dots.
These points represent the Pareto frontier obtained by a sampling average of the objectives at the Pareto optimal designs found by the methodology after the
fixed number of iterations, i.e., an estimation of $P[\mathbf{O}[\bx_{1:N_{\max}}]]$ which is to be contrasted to $P[\boldsymbol{\mu}_{N_{\max}}[\bx_{1:N_{\max}}]]$.
This Pareto frontier is a representation of the quality of the  solution obtained by the methodology.
With low levels of stochasticity the methodology neatly approximates the noise in the outputs as Gaussian, shown in \fref{ex1_a} (a) and (b).
With an increase in the value of the stochasticity parameter, $s$, the 
final Pareto frontier obtained starts diverging from $P_O$, shown in \fref{ex1_a} (c) and (d). 
In \fref{ex1_a} (c) and (d), the methodology ends up exploiting the area near 
the two ends of the observed $P[\mathbf{O}[\bx_{1:N_{\max}}]]$ only, and not the whole $P_O$ which is possibly a manifestation of the methodology not being able to estimate and filter out the excessive non-Gaussian noise in these cases.
The contours of the symmetric deviation (which can be understood as the probability of
a particular set of objective values being achievable conditional on the observations 
made thus far) do reinforce greater knowledge about the 
plausibility of the achievable values even in regions which tend to dominate the 
approximated Pareto frontier. This means that with
more simulations the methodology should eventually discover more 
Pareto efficient solutions across the complete boundary of the approximated Pareto frontier.
So, the symmetric deviation allows the decision maker to realize the potential 
value that lies in doing further simulations.

\subsection{Six-dimensional synthetic example}
\label{sec:validate_6d}
Consider the following test objective functions from \cite{knowles2006}:
\begin{gather}
    o_{1}(\bx,\omega) = \frac{1}{2}(x_{1}+s\xi_1(\omega))(1 + g), \\
    o_{2}(\bx,\omega) = \frac{1}{2}(1 - (x_{1}+s\xi_1(\omega)))(1 + g),\\
    g = 100\bigg[5 + \sum_{i \in \{2,\cdots,6\}}{((x_{i} + s\xi_{i}(\omega))-0.5)}^{2}\nonumber \\
    - \cos(2\pi((x_{i} + s\xi_{i}(\omega))-0.5))\bigg],
\label{eqn:dtlz1a}
\end{gather}
for $\bx \in X = [0,1]^6$, where $\xi_i\sim\mathcal{N}(0,1),i=1,\dots,6$ are independent.
As before, the expected objectives are not analytically available.
We use the same approximation technique as in the previous example to estimate the ground truth of $P_O$ for this test problem.
\begin{figure}[htbp]
    \subfigure[]{
        \includegraphics[width=0.50\textwidth]{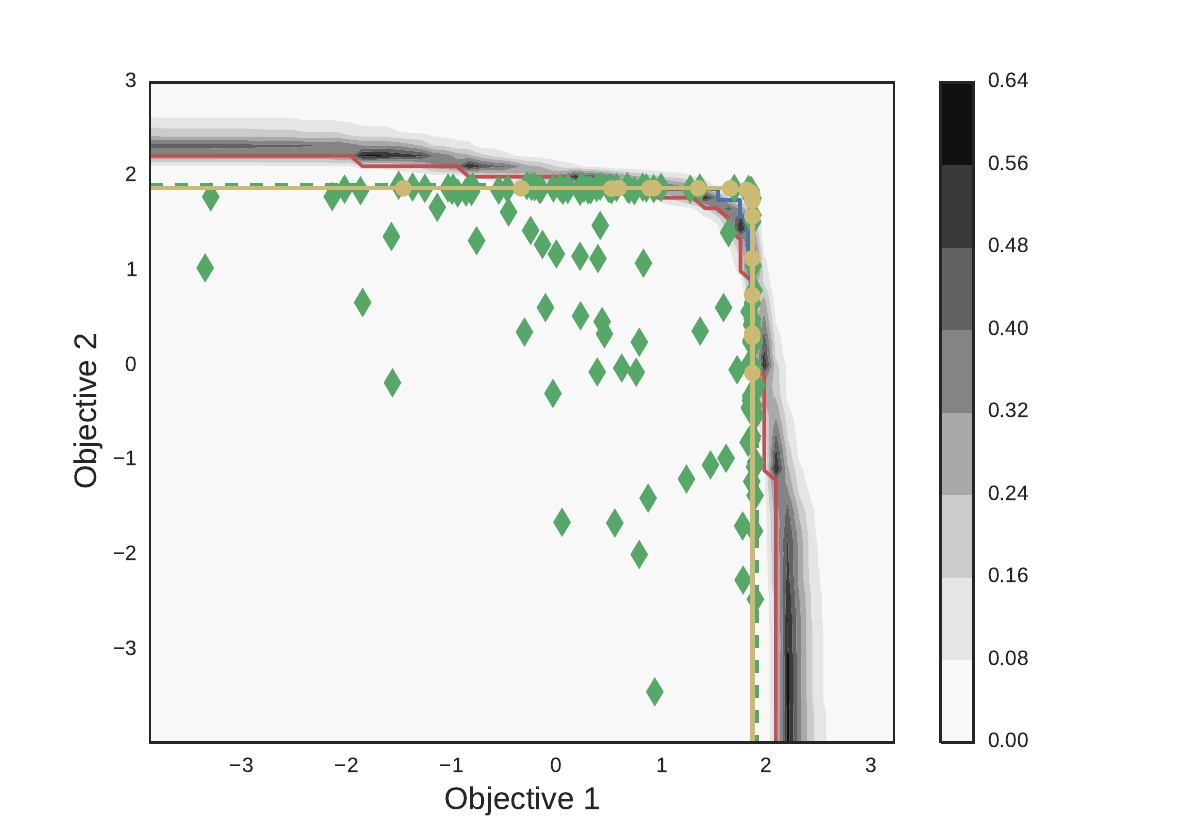}
    }
    \subfigure[]{
        \includegraphics[width=0.50\textwidth]{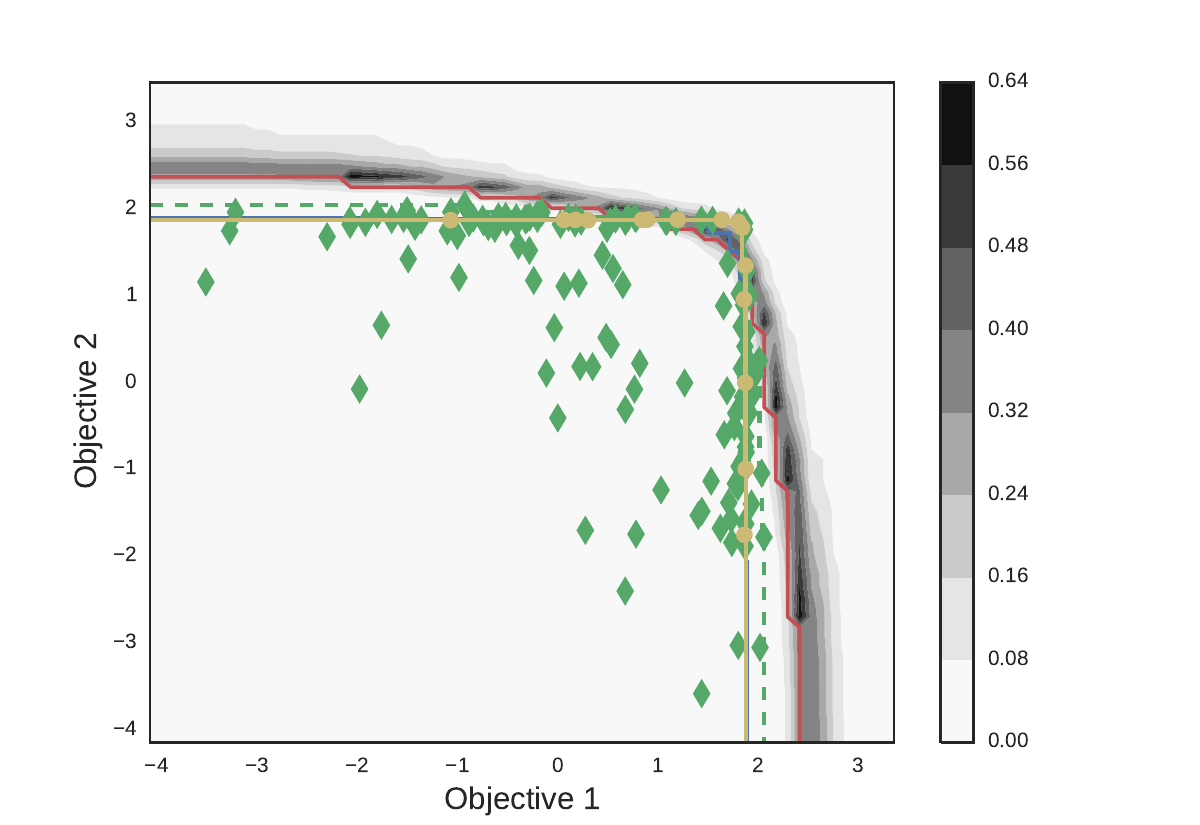}
    }
    \subfigure[]{
        \includegraphics[width=0.50\textwidth]{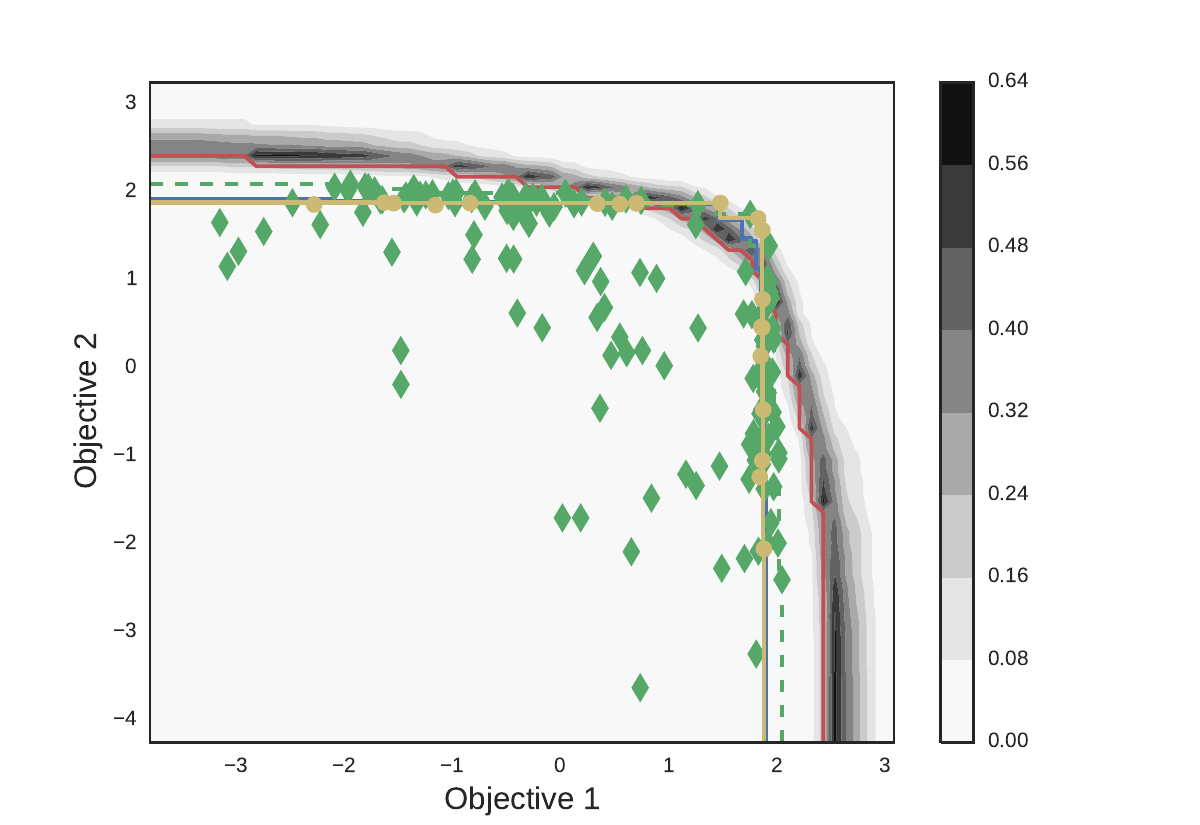}
    }
    \subfigure[]{
        \includegraphics[width=0.5\textwidth]{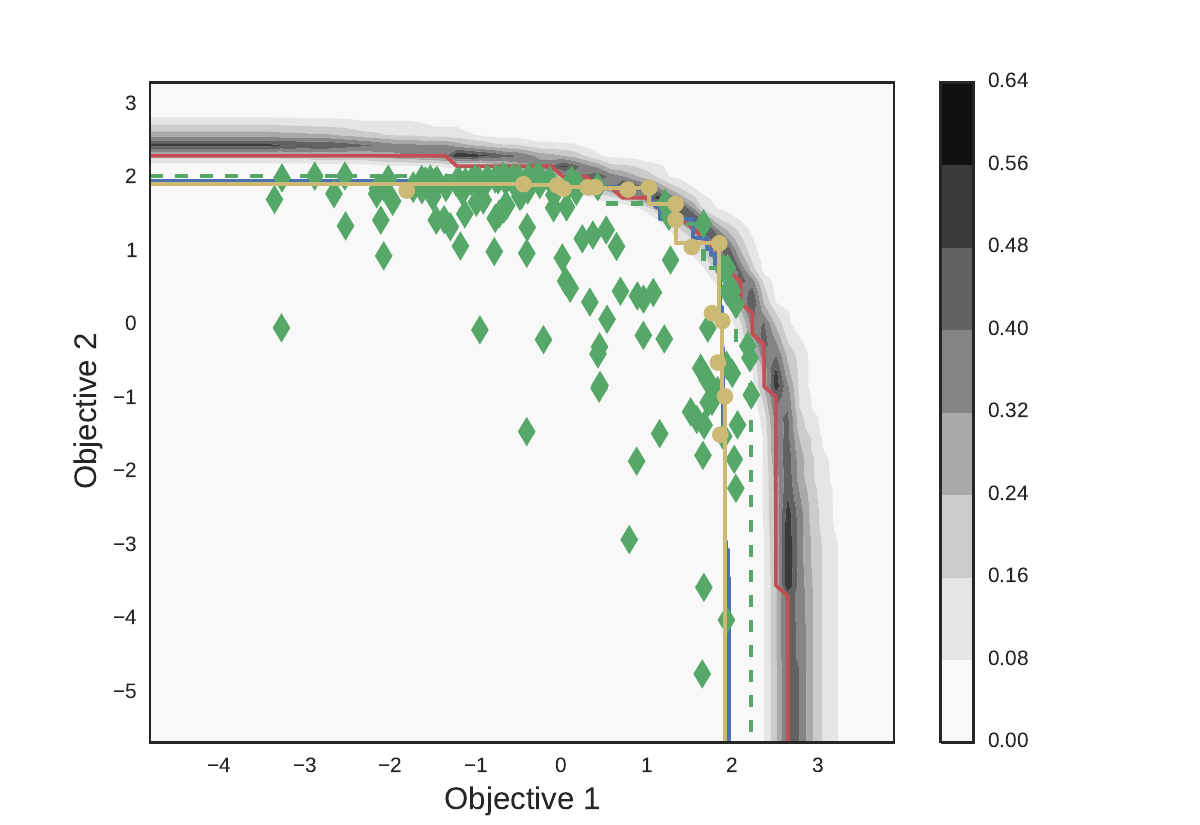}
    }
    \caption{Six-dimensional synthetic example starting from ($n=40$) initial measurements.
        Subfigures~(a)~($s=0.01$), ~(b)~($s=0.03$), ~(c)~($s=0.05$), and~(d)~($s=0.1$), 
        depict our state of knowledge about the final
        $P[\mathbf{f}^e[\bx_{1:n}]]$ after 100 measurements selected using Algorithm~\ref{alg:bgo}.
        }
    \label{fig:ex2_1a}
\end{figure}
\fref{ex2_1a} depicts our final state of knowledge about $P[\mathbf{f}^e[\bx_{1:n}]]$ for increasing noise levels $s = 0.01, 0.03, 0.05$, and $0.1$.
As before, the larger the noise the harder it is for the methodology to 
discover $P_O$, the true Pareto frontier.
In general, as can be seen in \fref{ex2_1a} the method is robust to noise as long as the noise is reasonably low for the given number of initial measurements. The powerfulness of the methodology can be observed through \fref{ex2_1a} (a) and (b) , where the final $P[\mathbf{f}^e[\bx_{1:n}]]$ contains points that dominate the  $P[\mathbf{O}[\bx_{1:N_{\max}}]]$, when the noise parameter has relatively low values.
The method, as expected, discovers very few points on $P_O$ as the noise increases to $s=0.1$ as can be seen in (d) of \fref{ex2_1a}.

%% file: wire-prob.tex
\subsection{Wire drawing problem}
\label{sec:wmp}
The wire drawing process is designed to achieve the desired final diameter and 
mechanical properties such as ultimate tensile strength (UTS) and ductility 
through cold reduction of a 
larger diameter wire. The desired wire properties depend on applications – 
for example, high torsional ductility is required for application in tires, 
high strength wires used in machine tools for metal cutting.
A typical reduction of the cross section the wire, based on the final properties 
required would be in the range of 70-90 percent
and this is achieved by reducing the wire diameter in a number of passes. 
Each pass involves drawing through a conical die and the sequence of reductions 
and corresponding die angles at each pass would play an important role on the 
final properties as well as performance of operations.
Here we consider a 
wire drawing process having a fixed number of passes (8 passes). An finite element analysis (FEA) based simulator, 
developed for an industrial operation was used to simulate this process. 
This wire drawing simulator includes wire deformation, heat generation and 
dissipation in the wire as well as dies, cooling of wire on the cooling drum 
and in the atmosphere and is based on large deformation theory.
The model considers the process to be axisymmetric.
The multi-pass drawing effect is modeled by considering carryover effect of previous pass such as residual stress, 
plastic strain and temperature.
The FEA is done using four noded isoparametric elements.
A penalty parameter approach is used for modeling the contact between the wire and the dies.
The simulator takes the input as wire material properties, input 
wire diameter, die pass schedule (reduction and die angle at each pass), 
wire drawing speed, cooling conditions, friction, etc.; and predicts the 
internal stress and strains in the wire and the die, load on each die and 
the drum, temperature of the wire and the die, properties indicative of final 
wire mechanical properties -- UTS representing strength and strain 
non-uniformity factor (SNUF) representing relative ductility.   

The plastic deformation across the cross section of the final wire should be as uniform as possible for enhanced ductility.
The UTS is primarily governed by the total reduction but the non-uniform deformation has a significant  secondary role on the final UTS.
To understand this uniformity, the plastic strain distribution is modeled and is represented as SNUF.
SNUF is a ratio of difference between the peak and average strain to average strain.
Besides the properties of the drawn wire, process defects such as wire burst during drawing process is an important aspect to consider as central burst is highly undesired since it leads to wire breakage during drawing process and this effect is modeled through the measurement of triaxiality by a factor called the hydraulic failure factor (HFF).
The coefficient of friction is assumed to be constant throughout the process.
Here, we have the UTS and the SNUF as the two competing objectives for the process. 

The design variables for this problem are the die angles (one at each pass) and the incoming wire diameter (implicit in the reduction ratio) at each pass. The outgoing wire diameter at a pass is same as the incoming wire diameter for the next pass.
The incoming wire diameter $d_j$ and the reduction ratio $(rr_j)$ for a pass $j$ are related by the formula given in \eqref{eqn:rr}.
\begin{equation}
    \label{eqn:rr}
    rr_j = 1 - \frac{d^2_{j+1}}{d^2_j}
\end{equation} 
For this problem we take the case of drawing an 8mm wire into a 3mm wire \fref{wire_draw}.
\begin{figure}[htbp]
    \centering
    \subfigure[]
    {
        \includegraphics[width=1.0\textwidth]{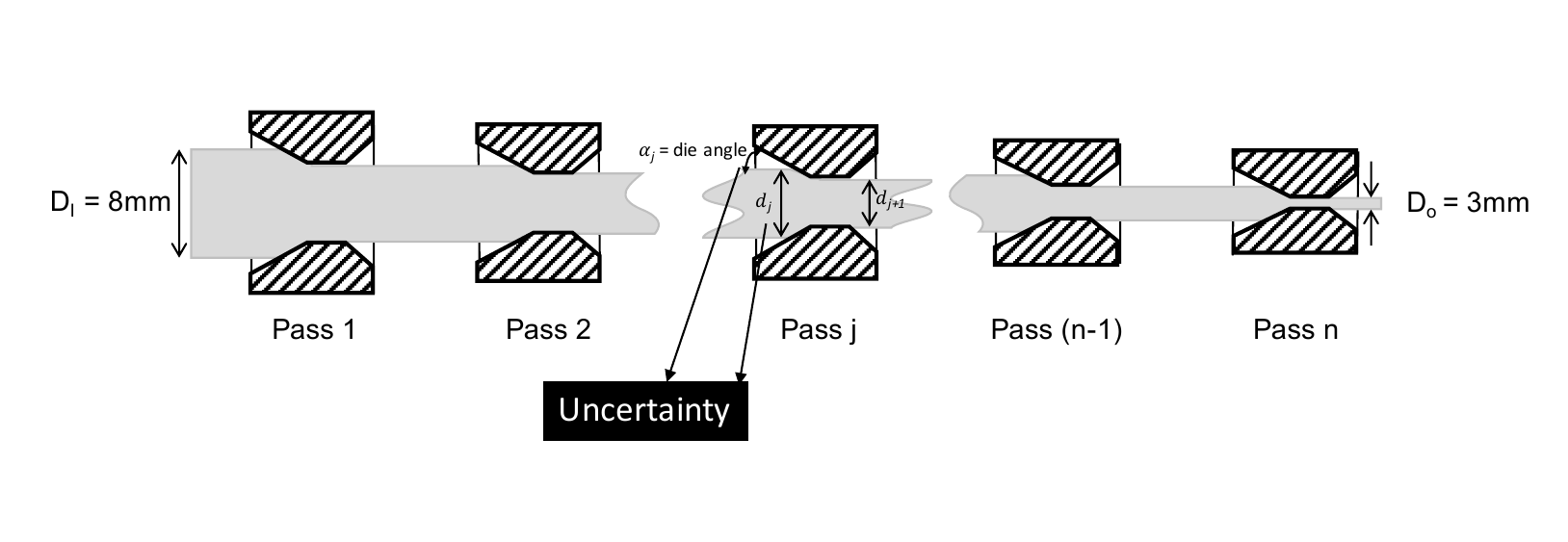}
    }
    \caption{WMP: The wire manufacturing process with the depiction of the sources of uncertainty, ie.
    the incoming wire diameter $d_j$ and the die angle $\alpha_j$, at an individual pass $j$. }
    \label{fig:wire_draw}
\end{figure}
So, with the overall reduction ratio (and the incoming wire diameter for the first pass) fixed, the problem becomes that of two objectives with 15 design parameters (8 die angles and 7 incoming wire diameters).
We apply our methodology to the wire drawing problem and demonstrate its ability to deal with the problem of stochasticity in the objectives induced by our inability to fully control the design parameters, to obtain a set of Pareto optimal solutions. This uncertainty can be understood as the ubiquitous effect of the continuous wear and tear on the die which would cause the process to deviate from delivering ideal (no noise) outputs. Also, in any manufacturing process the tolerances need to be accounted for as the procured dies themselves would not have exact dimensions as required.
The design space has been bounded by choosing a suitable range for design variables as follows:
\begin{enumerate}
    \item For $i=1,\dots,7$, $x_{i} \in [0,1] $ represent the \emph{incoming wire diameters}.
    \item For $i=1,\dots,8$, $x_{i+7} \in [0,1]$ represent the \emph{die angles}.
\end{enumerate}
Specifically, we assume that when we try to implement a process with design $\bx$, what we actually get is a process with design $\bx + \mathbf{S}\bxi$, where $\bxi\sim\mathcal{N}(\mathbf{0}_{15},\mathbf{I}_{15})$ and $S = \operatorname{diag}(s_1,\dots,s_{15})$ where $s_{i}=0.05, \forall i \in [1,7] $ and $s_{i}=0.1, \forall i \in [8,15] $. 
The above space $X = [0,1]^{15}$ is a scaled representation of the real space for simplification purposes.
The random vector from the real space $X = [7.2,7.5]\times[6.6,6.9]\times[5.8,6.1]\times[5.1,5.4]\times[4.4,4.7]\times[3.9,4.2]\times[3.3,3.6]\times[8,14]^8$, can be obtained by rescaling the random vector from the scaled space by using a simple linear transformation.
The noisy objectives considered here are:
\begin{eqnarray}
    o_{1}(\bx,\omega) &=& -\mbox{SNUF}\left(\bx + \mathbf{S}\bxi(\omega)\right),  \\
    o_{2}(\bx,\omega) &=& \mbox{UTF}\left(\bx + \mathbf{S}\bxi\right).
    \label{eqn:wmp-objectives}
\end{eqnarray}
The optimization problem involves maximizing the UTS and minimizing the SNUF.
For simplifying the problem to the requirements of our code and software
we convert it to an equivalent maximization problem where
we maximize the UTS and maximize the negative of the SNUF.
We consider a scenario with 15 initial observations of the MOO problem
and limit our computational budget to allow for 50 additional
simulations to be carried out sequentially. 
\begin{figure}[htbp]
    \centering
    \subfigure[]
    {
        \includegraphics[width=0.5\textwidth]{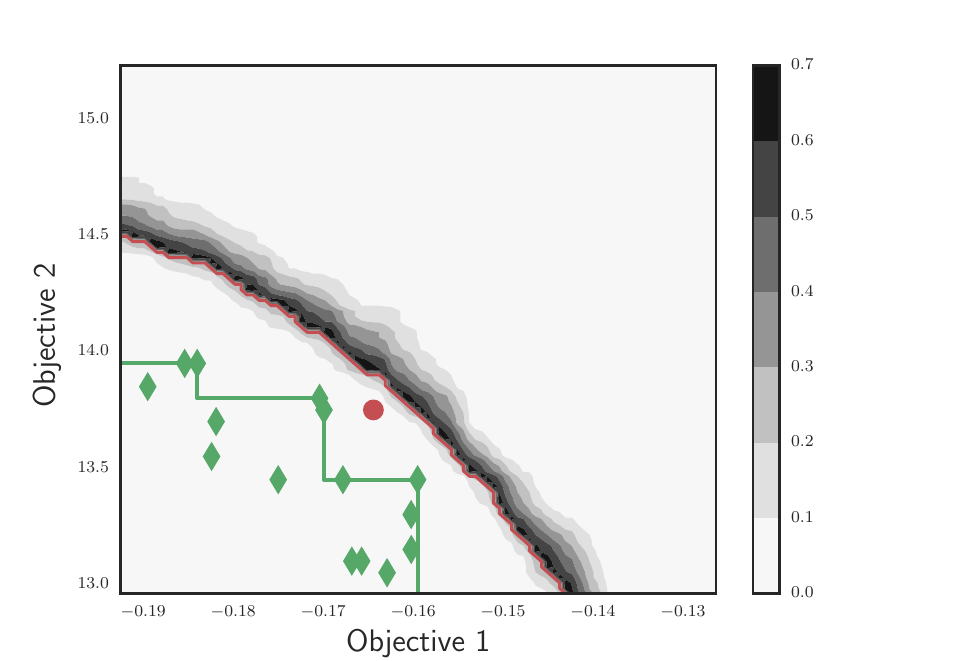}
    }
    \caption{WMP: The $P[\mathbf{f}^e[\bx_{1:n}]]$ for the inital observations using \qref{posterior_mean}. 
    Objective 1 is the -SNUF and Objective 2
        is the UTS.}
    \label{fig:wmp_1}
\end{figure}

\fref{wmp_1} shows the projected initial observations for the problem. 
We scale the measurements 
obtained by subtracting and dividing by the empirical mean and standard deviation just as in the case of the
test function discussed above.
This is done to maintain consistency with the assumption (in Sec. \ref{sec:prior}) of a zero mean (standard normal) 
GP for computational flexibility.

A key aspect of quantifying our knowledge about the state of the objectives is the Vorob'ev expectation which
is computed by obtaining by sampling the design space $X$. However, it must be noted that in this case with 15 dimensions,
it becomes very difficult to cover the whole design space as a result of which certain designs picked by the algorithm, end up outside the sampled designs. The overarching effect of this can be seen in \fref{wmp_2} (a), where the Vorob'ev expectation can be seen lying below the points in the top left corner picked by the methodology. To circumvent this issue, we augment the set of  
sampled designs with the designs at which we have made observations. This provides a clearer picture, \fref{wmp_2} (b), of the state as it reinforces the information obtained thus far while quantifying our beliefs about the state of the Pareto-efficient frontier.  
\begin{figure}[h!]
    \subfigure[]{
        \includegraphics[width=0.5\textwidth]{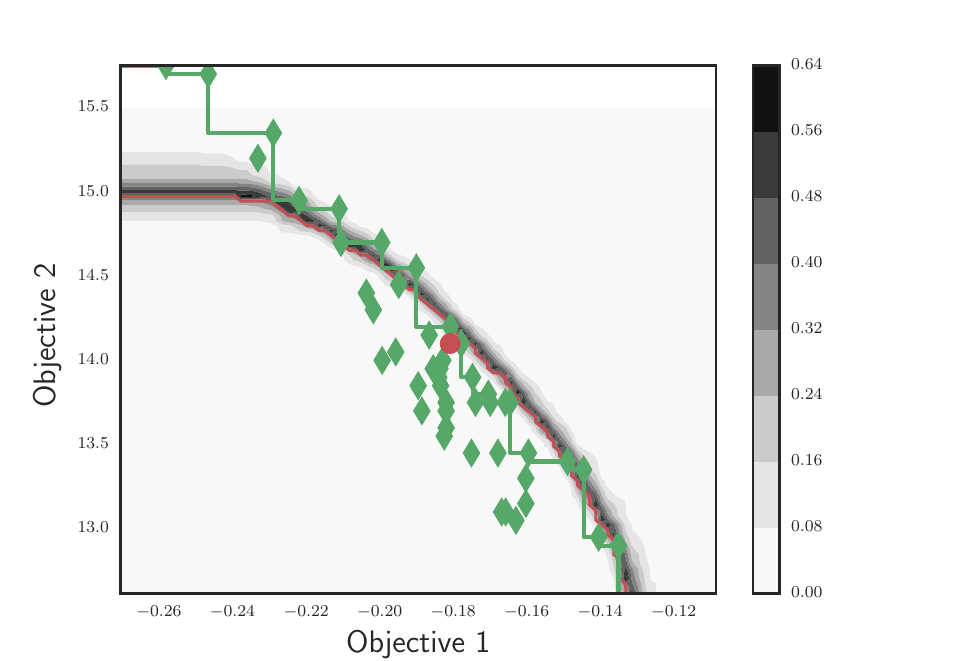}
    }
    \subfigure[]{
        \includegraphics[width=0.5\textwidth]{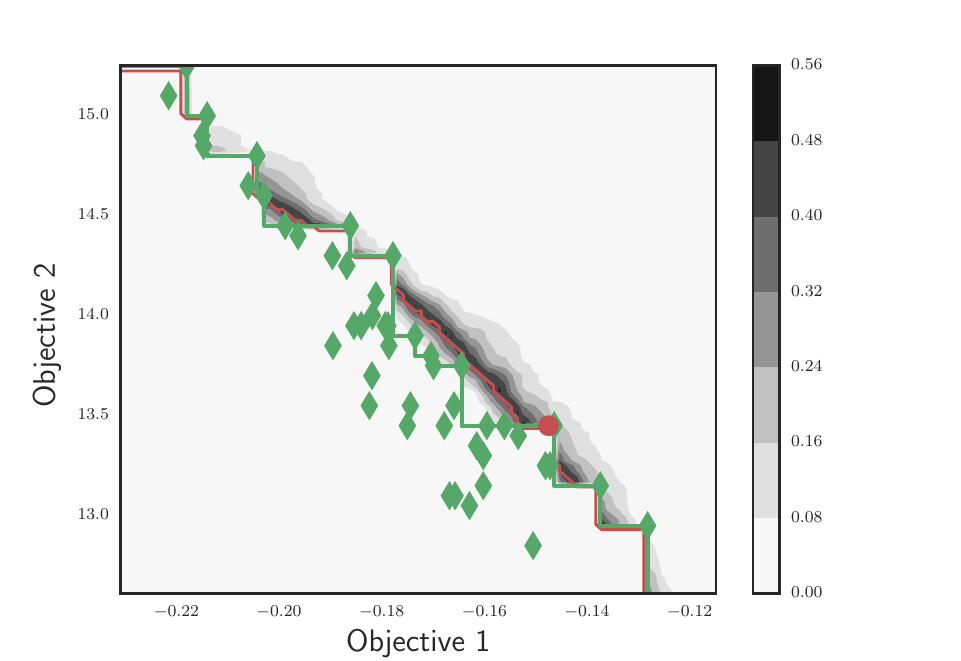}
    }
    \caption{WMP: The $P[\mathbf{f}^e[\bx_{1:n}]]$. Subfigures~(a)~ with the random sample design space, ~(b) after adding the observed designs to the sampled design space.
        }
    \label{fig:wmp_2}
\end{figure}

\fref{wmp_3} depicts the state of the
problem after the fiftieth iteration. Since, we do not have the computational
resources to obtain  $P[\mathbf{O}[\bx_{1:N_{\max}}]]$ for comparison, we sample average the value of the
objectives, 100 times, corresponding to the final Pareto designs as shown in \fref{wmp_2}. This averaging gives us an
estimate of the approximate true state of the Pareto-efficient frontier after the computational budget has been exhausted.
\begin{figure}[h!]
    \centering
    \subfigure[]
    {
        \includegraphics[width=0.5\textwidth]{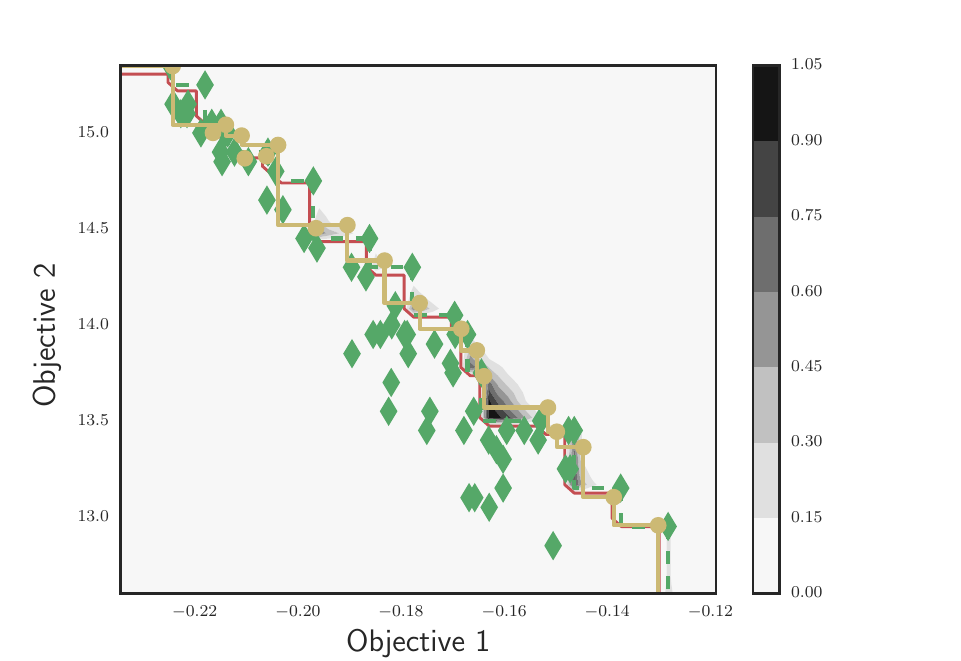}
    }
    \caption{WMP: The $P[\mathbf{f}^e[\bx_{1:n}]]$ after 50 additional measurements
        along with the \emph{sampled averaged approximation of $P[\mathbf{O}[\bx_{1:N_{\max}}]]$} represented by the yellow line. 
        Objective 1 is the -SNUF and Objective 2
        is the UTS.}
    \label{fig:wmp_3}
\end{figure}

%% file: conclusions.tex
\section{Conclusions}
\label{sec:conclusions}
{We constructed an extension to the EIHV information acquisition function which
makes possible the application of BGO to stochastic
multi-objective black-box optimization problems.
In addition to the above, we have shown how the epistemic uncertainty induced by the limited number of simulations can be quantified and 
used, to represent the uncertainty around the PF at each stage. 
We have validated our approach by applying it on two, slightly modified to include stochastic parameters, synthetic test functions with known Pareto frontiers.
Furthermore, we applied our method on the challenging steel wire drawing problem under parametric uncertainty in a scenario of simulation based design. 
The method offers a viable alternative to the state-of-the-art evolutionary optimization algorithms which rely heavily on sample averaging and are unaffordable
under a limited budget scenario.
Moreover, the proposed extension to EIHV gives acceptable results under cases of moderate levels of noise with limited number of initial observations.
There remain several open research questions. 
The most pressing direction to look in would be the efficient treatment of stochastic multi-objective problems under unknown and expensive constraints under a scenario of constrained computational resources.}